\documentclass[twocolumn]{article}
\usepackage[utf8]{inputenc}
\usepackage{cite}
\usepackage{amsmath, amsfonts, amssymb, amsthm}
\usepackage{graphicx}
\usepackage[frozencache=true, cachedir=minted-dir]{minted}
\usepackage{geometry}
\usepackage{csquotes}
\usepackage{array}

\usepackage{thm-restate}
\newtheorem*{ack}{Acknowledgements}

\newcommand{\Sim}{\textbf{Sim}}

\title{\textbf{Automatic differentiation as an effective tool in Electrical Impedance Tomography}}
\author{\textbf{Ivan Pombo$^{1, 2}$, Luis Sarmento$^2$} \\  	\\
	$^1$ CIDMA - Center for Research and Development in Mathemathics and Applications 
	\\ $^2$ Inductiva Research Labs \\
	\{ivanpombo.ext, sarmento\}{\fontfamily{ptm}\selectfont @}inductiva.ai
} 
\date{ }

\begin{document}
	
\maketitle

\begin{abstract}
    Determining physical properties inside an object without access to direct measurements of target regions can be formulated as a specific type of \textit{inverse problem}. One of such problems is applied in \textit{Electrical Impedance Tomography} (EIT). 
    
    In general, EIT can be posed as a minimization problem and solved by iterative methods, which require knowledge of derivatives of the objective function. In practice, this can be challenging because analytical closed-form solutions for them are hard to derive and implement efficiently.
    
    In this paper, we study the effectiveness of \textit{automatic differentiation (AD)} to solve EIT in a minimization framework. We devise a case study where we compare solutions of the inverse problem obtained with AD methods and with the manually-derived formulation of the derivative against the true solution.
    
    Furthermore, we study the viability of AD for large scale inverse problems by checking the memory and load requirements of AD as the resolution of the model increases.
    With powerful infrastructure, AD can pave the way for faster and simpler inverse solvers and provide better results than classical methods.
\end{abstract}

\section{Introduction}

Electrical Impedance Tomography (EIT) is a non-invasive
imaging method that produces images by
determining the electrical conductivity inside a
subject using only electrical measurements obtained at its
surface. More specifically, sinusoidal currents are applied
to the subject through electrodes placed in certain locations at the surface of the object. The
resulting voltages are then measured, making it possible to infer internal properties of the objects. EIT is a low-cost method and harmless for human being, since it only applies low amplitude currents. Additionally, it allows for real-time monitoring of various subjects even in the most difficult conditions. There are
applications of this technology for medical purposes, in scenarios such as ventilation monitoring, detecting brain hemorrhages and breast cancer. EIT is also used in geophysical imaging, flow analysis and other industrial purposes. For further insight into the applications see \cite{webster} and \cite{adler}. 

\begin{figure}[ht]
	\centering
	\includegraphics[scale=0.6]{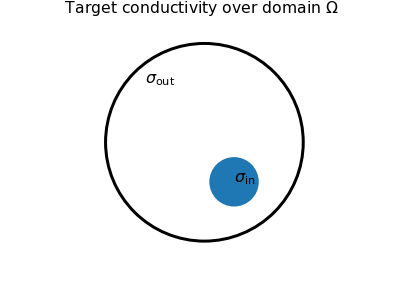}
	\caption{Example of a target conductivity over the domain $\Omega$ that represents a simple model of breast cancer where tumors have higher conductivity than the background. The domain $\Omega$ is represented by the black circumference which has a conductivity of $\sigma_{\mathrm{out}}$. In a blue circle it is represented a region with different conductivity $\sigma_{\mathrm{in}}$ from the background one $\sigma_{\mathrm{out}}$. }
	\label{fig:cancer_id}
\end{figure}

A particularly relevant application of EIT is in the early determination of breast cancer, specifically for young women where the risk of the ionizing X-rays of mammographies outweigh the benefits of regular check-ups. Fig. \ref{fig:cancer_id} describes one simplified EIT scenario where the blue region represents cancer inside the breast, denoted as a domain $\Omega$. The assumption is that healthy and cancerous tissue have different conductivity values $\sigma_1, \, \sigma_2$, respectively. The goal is to locate a potential region affected by cancer from measurements on the breast surface, which is the boundary of the domain $\Omega$ and denoted as $\partial\Omega$.

The measurements are obtained by injecting into the domain $\Omega$ a fixed set of different \textit{electrical current patterns} $I_j$. Each $I_j$ is defined by injecting electrical current through all electrodes in a particular manner, \textit{i.e.}, for $L$ electrodes we have $I_j=(I_{j,1},...,I_{j,L})$. Simultaneously, we measure the resulting voltages $V_j$ for each current pattern,  obtaining a voltage measurement at each electrode, denoted as $V_j=(V_{j,1},...,V_{j,L})$. This leads to a set of \textit{true measurements} denoted by $m_j=(I_j, V_j)$. Then, the corresponding \textit{inverse problem} is to determine the electrical conductivity over $\Omega$ that leads to these measurements. In the particular case of Fig. \ref{fig:cancer_id} we want to determine the conductivity outside and inside the anomaly, $\sigma_{\mathrm{out}}$ and $\sigma_{\mathrm{in}}$, respectively, and the location of the anomaly (in blue).

This is a hard problem because in general there is no analytical expression that maps a set of electrical measurements back to the respective conductivity profile that generates them. 

To solve this inverse problem we first need to understand how to solve the \textit{direct problem}, that is, computing electrical measurements $V_j$ for a given set of currents $I_j$ and conductivity $\sigma$. The direct problem has an easier solution, since the propagation of electrical current through the domain obeys the well-known Maxwell equations. 

Many methods for solving the direct problem are described in the literature, \textit{e.g.}, Finite Element Method (FEM) \cite{vauhkonen}, Boundary Element Method (BEM) \cite{tanzer}, and, more recently Deep Learning methods (DL) \cite{raissi}.

Independently of the numerical method used to solve the direct problem, such a procedure is commonly designated as \textit{simulation}. Hence, for a given conductivity profile we can obtain through a simulation method the electrical measurements denoted as $m^{\text{Sim}}_j=(I_j,V^{\text{Sim}}_j), $ for each different current pattern with $j=1,...,N$. We can thus define an operator that maps conductivity into voltage measurements, here termed by \textit{direct operator}and given as:
\begin{align}\label{sim}
	\Sim: \sigma\mapsto  V^{\Sim}=(V_{1,1}^{\text{Sim}},.., V_{j,l}^{\text{Sim}},...,V_{N,L}^{\text{Sim}})\in\mathbb{R}^{L\cdot N}
\end{align}
where $V_{j,l}^{\text{Sim}}$ represent voltages measured at the $l$-th electrode for the $j$-th current pattern.

Our goal is to find a conductivity profile $\sigma$ that matches measurements $m=(m_1,...,m_N)$. Thus, we can formulate EIT as the following minimization problem by making use of the direct operator $\Sim$:
\begin{align}\label{min}
	\min_{\sigma} \frac{1}{2}\left\|\Sim(\sigma) - m^{\text{true}}\right\|_2^2.
\end{align}

We use the $L^2$-norm here for simplicity, but, in general, we could use any other norm as long as it is differentiable.

Most classical methods for solving this minimization problem are based on iteratively improving the solution. The update requires computing the derivative of both the loss function in (\ref{min}) and the $\Sim$ operator. 

To solve the inverse problem under an optimization framework we opted for the \textit{Levenberg-Marquardt algorithm} \cite{levenberg, marquardt}. It is a simple quasi-Newton method that only requires the Jacobian computation of the $\Sim$ operator. Further details about the method are given in Appendix B.

In essence, the main challenges to solve the minimization problem (\ref{min}) with iterative classic methods are: 
\begin{itemize}
	\item to ensure that the simulator is once-differentiable with respect to a conductivity parameterization; 
	\item devise a method to compute the respective derivatives of the simulator.
\end{itemize}

Our study explores a simulation operator obtained through FEM, which is already well established for EIT, see \cite{siltanen}. 
\begin{figure}[h]
	\centering
	\includegraphics[scale=0.6]{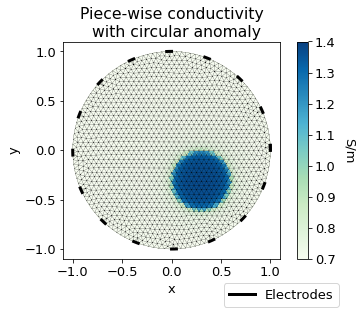}
	\caption{Circular anomaly defined over a triangular FEM mesh in 2D. Electrodes are attached to the boundary, black lines.}
	\label{fig:mesh_anomaly}
\end{figure}

When the $\Sim$ operator is given by FEM we can deduce an analytical closed-form of the derivative with respect to the conductivity variation. It is simply obtained with respect to a conductivity discretization over the FEM mesh, see Fig. \ref{fig:mesh_anomaly}. As such, it requires derivative computations with respect to conductivity values over \textit{all elements of the mesh}. If the conductivity is defined through a different parameterization we can obtain the respective derivatives by the chain rule of differentiation. For such endeavor, the analytical formulation needs to be adapted and derived for each particular parameterization of the conductivity. As a result this formula is hard to derive and implement, see \cite{harrach} and \cite{other_vauhkonen}.

Automatic differentiation (AD) is a method that automatically evaluates exact derivatives for complex programs. It exploits the simple mathematical operations the programs are built on, to automatically compute the derivative through the chain rule. While the initial concept was developed in the sixties \cite{wengert}, only recently with advancements in hardware and efficient implementations, like JAX \cite{jax}, it has gained traction for application in general problems. 

In this paper, we explore automatic differentiation as an alternative to manual methods for computing the Jacobian of differentiable simulators. In particular, the goal is to validate its effectiveness in solving the EIT inverse problem. By effectiveness we mean that it is as successful in solving the inverse problem as previous methods, namely, through analytical formulation. By doing so we show its versatility compared with analytic formulation and moreover verify its viability for high resolution images where . 

The validation is done by comparing the absolute error between solutions obtained by solving the minimization problem with both methods to compute the derivatives and the absolute error compared with the true solution. In particular, we evaluate the maximum difference between both Jacobian computations to check if they are evaluating to the same result. Then as a second set of checks, we explore the memory consumption of AD and show that it is still in reasonable terms as the problem scales with higher resolution.

Our end goal is to show feasibility and practicality of AD as a tool for lowering the entry barrier for other inverse problems in Partial Differential Equations, where AD can also be applied.

In the following, we first introduce the EIT case study we are using for comparison. In Section 3, we explain how the required derivatives are computed with both the analytical and AD method. In Sec. 4, we introduce our experimental setup. Results comparing the effectiveness of both methods and viability of AD are given in Sec. 5, and conclusions are drawn out in Sec. 6. 

\section{Establishing a case study} 

In this section we establish a case study in order to make a clear comparison between both methods for computing derivatives.

\subsection{EIT scenario}

To demonstrate our claims we focus on a two-dimensional setup. We remark that this is not physically accurate since electrical current propagates in three dimensions However, it simplifies the construction of our case study.

EIT is an \textit{ill-posed} inverse problem \cite{siltanen} and thus we need to take into account the possible instability of the problem, \textit{i.e.}, small variations in the measurements may imply large variation on parameters solution. In practice, this makes it hard to solve the inverse problem since true measurements, captured with real-world measuring devices, always contain noise. Therefore, real solutions for noisy input data can be drastically distinct from the true solution.

Due to this, it becomes hard to accurately determine a very large number of parameters, \textit{e.g.}, the value of the conductivity at all mesh elements (see Fig. \ref{fig:mesh_anomaly}), from a small number of measurements. An example would be a conductivity defined over a fine mesh which has a value at each mesh element, see Fig. \ref{fig:mesh_anomaly}. 

To mitigate this problem we want to make as many measurements as possible. However, the possible number of distinct measurements is constrained by the quantity of electrodes. This occurs since for $L$ electrodes there are only $L-1$ linearly independent current patterns for which the voltage measurements yield independent information of the conductivity, see \cite{siltanen}. 

The best way to mitigate this issue is to work on simpler cases. By doing so we can reduce the parameter space and have less variability on the solutions,  like in Fig. \ref{fig:cancer_id}.
Even though instability issues do not become completely fixed, the space of measurements has a lower, more tractable, dimensionality. 

For the sake of comparison we wish to make, it is enough to focus on conductivity profiles with a circular region of distinct conductivity value from the background, see Fig. \ref{fig:cancer_id} and \ref{fig:mesh_anomaly}. These anomalies are parameterized by their center $(c_x, c_y)$ inside the domain $\Omega$, radius $r$ and conductivity value inside and outside $\sigma_{\mathrm{in}},\, \sigma_{\mathrm{out}}$, respectively.

We work with this simplification because  it is easier to obtain a solution to the inverse problem due to the parameterization of such region being given by only a few parameters. Further, we remark that we need to make sure that the parameterization is differentiable. Our choice of circular regions is based on this, since it is easy to define a smooth parameterization. For regions with corners two smoothing procedures would be required, one to smoothen the corners and another to smooth the parameterization.

By the reasons above, in our experiments we assume the existence of a \textit{single} circular anomaly with conductivity value different from the background, like in Fig. \ref{fig:cancer_id}and denoted the parameterization variables as
\begin{align}\label{anomaly}
	\sigma = (r, c_x, c_y, \sigma_{\mathrm{in}}, \sigma_{\mathrm{out}}).
\end{align}

 We introduce now the EIT model, the conductivity parameterization definition and the measurement setup we use to proceed with out comparison.

\subsection{Voltage measuring setup} 

We introduce here the measuring setup that is applied for the direct problem.

In this simple 2D setup, we define the $\Sim$ operator in (\ref{sim}) according to the case study and the measurement setup.

Recall that with $L$ electrodes at the surface $\partial\Omega$, we can at most apply $L-1$ linearly independent current patterns $I_j\in\mathbb{R}^L$ with $j=1,...,L-1$. The $\Sim$ operator is obtained by solving the direct problem for each $I_j$ and determine the respective voltages $V_j\in\mathbb{R}^L$ over the electrodes. 

The more measurements we can perform the better we are able to potentially reconstruct the conductivity. Therefore, we need to choose $L-1$ linearly independent current patterns. This choice is non-trivial. One possibility presented in the literature \cite{siltanen} is obtained by injecting currents in a wave pattern through the electrodes according to
\begin{align}\label{current_pattern}
I_{j, l} = \begin{cases}
A \cos(j\theta_l),\;\; &j=1,..., \frac{L}{2},\\
A\sin\left((j-\frac{L}{2})\theta_l\right),\;\; &j=\frac{L}{2}+1,...,L-1
\end{cases}
\end{align}
with $\theta_l = \frac{2\pi}{L}l$ and $A$ the constant current amplitude. These patterns have been shown to obtain th best result on the detection of conductivities profiles with small anomalies in the regions furthest from the boundary, \cite{siltanen}.

The experiments are performed in the following setting:
\begin{itemize}
	\item $\Omega$ is a circular domain with radius  $r_{\Omega}10\mathrm{cm}$;
	\item Current amplitude of $A=3\mathrm{mA}$, which is a reasonable value for human subjects, and the voltages are measured in (mV);
	\item Attach $L=16$ electrodes equally spaced at the boundary with each having fixed length $\pi/64$. 
\end{itemize}

We refer to Figure \ref{fig:mesh_anomaly} for a visual representation of the setting.

Under the above setup, the simulator in equation (\ref{sim}) is given as
\begin{align}\label{para_Sim}
	\Sim: \mathbb{R}^5 &\rightarrow \mathbb{R}^{L(L-1)} \\ \nonumber (r, c_x,c_y, \sigma_{in}, \sigma_{out})&\mapsto (V^{\Sim}_1,...,V^{\Sim}_{j,l},...,V^{\Sim}_{L-1, L})
\end{align} 
with $V^{\Sim}_{j,l}\in \mathbb{R}$ being the voltage measurement on the $l$-th electrode obtained by the direct problem solution for the trigonometric current pattern $I_j$.

\section{Modeling EIT}

\subsection{Direct problem}

Currents propagating in human tissues and organs can be satisfactorily modeled by the Complete Electrode Model \cite{cheng}. It accounts for the finite nature of electrodes, for the current injection through them and for the electro-chemical effects happening between skin and electrode surface. 

Let $\Omega$ describe the subject region we are evaluating. To establish a measurement setup, we attach $L$ electrodes at the subject boundary $\partial\Omega$. Through them we apply an electrical \textit{current pattern} $I=\left(I_1,...,I_L\right)$ into $\Omega$. The objective is to find the electrical potential $u$ inside and the voltages at electrodes $V=\left(V_1,...,V_L\right)$
that fulfill the system of equations describing the Complete Electrode Model:
\begin{align}\label{CEM}
\begin{cases}
\nabla\cdot(\sigma\nabla u)=0, &\quad\text{ in } \Omega,\\
\int_{E_l} \sigma\frac{\partial u}{\partial \nu}\, dS = I_l, &\quad l=1,2,...,L \\
\sigma\frac{\partial u}{\partial\nu} = 0,&\quad \text{ in } \partial\Omega\setminus\cup_{l=1}^L E_l \\
\left. u + z_l\sigma\frac{\partial u}{\partial\nu}\right|_{E_l} = V_l, &\quad l=1,2,...,L
\end{cases}
\end{align}
where $\nu$ is the outward pointing normal vector at $\partial\Omega$, $dS$ is measuring length of the boundary and $\sigma$ is the conductivity distribution. 

The first equation represents electrical \textit{current diffusion}. The second and third define the insertion of current through electrodes, meaning current spreads through the whole electrode before being inserted into the domain and in regions without electrodes there isn't current flowing. Finally, the last equations model the electrochemical effects at interface of skin-electrode, with $z_l$ termed as \textit{contact impedance} representing the resistance at that interface.

To ensure the existence and uniqueness of a solution, the current pattern must satisfy Kirchoff's law and we fix a \textit{reference voltage condition}:
\begin{align}\label{uniq_cond}
\sum_{l=1}^L I_l=0, \quad \text{ and }\quad \sum_{l=1}^L V_l=0.
\end{align}

\subsection{Modeling the circular anomaly}

In this section, we define the conductivity parameterization formally introduced in Section 2.

The parameterization is done through a \textit{level-set}, \textit{i.e.}, a function that has positive sign inside the region it describes, negative on the outside and equal to zero on the region boundary. In particular, a \textit{circle level-set} LS$(x,y)$ can be defined through a center $c=(c_x,c_y)$ and a radius $r$ as follows
\begin{align}
\text{LS}(x,y)=r^2 - \left[(x-c_x)^2+(y-c_y)^2\right].
\end{align}

The level-set function is positively valued if the point $(x,y)$ is inside the circular anomaly, negative if it is outside and zero if its precisely at the boundary of the anomaly.

As such, we can use the \textit{Heaviside function} $H(z)$ that equals $1$ if $z>0$ and $0$ otherwise, to fully describe the conductivity profile of interest through
\begin{align}\label{conductivity}
\sigma(x, y) = \sigma_{in} H(\text{LS}(x,y))  + \sigma_{out}\left(1-H(\text{LS}(x,y))\right).
\end{align}

Under this formulation $\sigma$ is not differentiable due to the discontinuity of $H$ at $z=0$. In order to attain differentiability, we use a smooth approximation of the Heaviside function given as
$$H^{\epsilon}(z) = \frac{1}{\pi} \arctan\left(\frac{z}{\epsilon}\right)+\frac{1}{2}.$$
The conductivity $\sigma$ is instead established in terms of $H^{\epsilon}$, where $\epsilon>0$ works as a smoothing parameter. The smaller it is the closer $H^{\epsilon}$ is to $H$. 

This smoothing procedure is necessary both for the analytical computation as well as AD. In fact, we need to take into account the mathematical differentiability for a proper implementation of derivatives through AD. For example, JAX AD applies the derivative to $H$ by following the conditional operations \texttt{if else}, which implies a derivative of $0$ everywhere, which is not true for $z=0$.

\section{Derivatives computation}

In order to solve the inverse problem in a minimization framework, we need to compute derivatives of the $\Sim$ operator. In this section, we deduce the analytical formula and explain how to apply AD to $\Sim$, in order to obtain the derivatives with respect to the parameters of interest. 

We recall that the direct solver and $\Sim$ are independent of the derivative computation method.

\subsection{Analytical Computation}

We recall that by Eq. (\ref{para_Sim}) we have that the FEM simulator operator is given by
\begin{align}\label{fem_sim}
	\Sim: \mathbb{R}^5 &\rightarrow \mathbb{R}^{L(L-1)} \nonumber \\
	(c_x,c_y, r,\sigma_{in}, \sigma_{out})&\mapsto \left(V^{\Sim}_1,...,V^{\Sim}_{j,l},..., V^{\Sim}_{L-1,L}\right).
\end{align}

To avoid heavy notation, we denote the vector of voltage measurements by $V^{\Sim}\in\mathbb{R}^{L(L-1)}$ and $V_{n}\in\mathbb{R}^L$ are the voltages measured $j$-th current pattern.

The Jacobian matrix $J\in \mathbb{R}^{L(L-1)\times 5}$ is given by
\begin{align}\label{jacobian}
	J=\begin{pmatrix}
	\frac{\partial V^{\Sim}}{\partial c_x} & \frac{\partial V^{\Sim}}{\partial c_y} & \frac{\partial V^{\Sim}}{\partial r} & \frac{\partial V}{\partial \sigma_{in}} & \frac{\partial V^{\Sim}}{\partial \sigma_{out}}\\
	\end{pmatrix} 
\end{align}

In order to provide an analytical formulation, we specifically focus on the computation of derivatives for each $V_n$ with respect to a single parameter, which if done for all $n=1,...,L-1$ determines one column of the Jacobian.

Furthermore, we need to specify a method to simulate the measurements. 

In this paper, we have used FEM applied to the Complete Electrode Model described before. The FEM solution is $\theta = (\alpha, \beta) \in \mathbb{R}^{N+L-1}$, where $\alpha$ describes the electrical potential inside $\Omega$ and $\beta$ the voltages at the electrodes. Accordingly, we denote for each current pattern $I_j$ the FEM solution by $\theta_j=[\alpha_j, \beta_j]\in\mathbb{R}^{N+L-1}$ with respect to $\tilde{I}_j$ on the right-hand side of the FEM system of equations (a variation of $I_j$). 

With this in mind, the voltages are computed by $V_j=M\beta_j$ where $M$ is a matrix defining the basis functions used by FEM at the electrodes. For further detail about the FEM solution we point to Appendix A.

Now, if we define $\tilde{M}=[\hat{0} \; M]\in\mathbb{R}^{L\times(N+L-1)}$ then we have 
\begin{align}
V_n=\tilde{M}\theta_n = \tilde{M}A^{-1}\tilde{I}_n.
\end{align}

As such, it holds for any parameter $w$ of $\{c_x,c_y,r,\sigma_{in},\sigma_{out}\}$ that:
$$\frac{\partial V_n}{\partial w} = \frac{\partial\left(\tilde{M}A^{-1}\tilde{I}_n\right)}{\partial w}.$$
	
Since neither $\tilde{M}$ and $\tilde{I}_n$ depend on the conductivity $\sigma$ and, therefore, for any of the parameters, it holds that
\begin{align}\label{der_analytic}
	\frac{\partial V_n}{\partial w}=\tilde{M}\frac{\partial A^{-1}}{\partial w}\tilde{I}_n = - \tilde{M}A^{-1}\frac{\partial A}{\partial w}A^{-1}\tilde{I}_n
\end{align}
with the last equality following from matrix calculus properties.

Thus, in essence, the computation resumes to the stiffness matrix derivative and noticing that $A^{-1}\tilde{I}_n=\theta_n$. Setting $\gamma=\tilde{M}A^{-1}$ the computation of the derivative in Eq. (\ref{der_analytic}) simplifies to
\begin{align}\label{analytic_der}
	\frac{\partial V_n}{\partial w} = -\gamma^T \frac{\partial A}{\partial w}\theta_n.
\end{align}

As such, the focus is on the computation of $\frac{\partial A}{\partial w}$. The stiffness matrix $A$ is composed of four blocks, like,
\begin{align*}
	\begin{bmatrix}
		B^1 + B^2 & C \\ C^T & D
	\end{bmatrix}.
\end{align*} 
The block $B^1$ is the only one depending on the conductivity. Due to its definition there is a clear way of computing the derivatives of $B^1$ with respect to the conductivity value  $\sigma_k$ over each mesh element (see the Appendix for further details on its definition):
\begin{align}\label{b1_derivative_sigma_k}
	\frac{\partial B^1_{ij}}{\partial \sigma_k} = \begin{cases}
	\int_{T_k} \nabla\phi_i\cdot\nabla \phi_j\, dx, \text{ if } i,j \in T_k \\
	0, \text{ otherwise}.
	\end{cases}
\end{align}  

Furthermore, the resulting matrix is independent of $\sigma$ therefore it can be precomputed at the start and re-used.  

Through the chain rule we have that
\begin{align}\label{b1_der_w}
	\frac{\partial B^1_{ij}}{\partial w} = \sum_{k=0}^K \frac{\partial B^1_{ij}}{\partial \sigma_k}\frac{\partial \sigma_k}{\partial w}.
\end{align}

We note that due to sparsity of the matrix defined in Eq. (\ref{b1_derivative_sigma_k}) it can be assembled very efficiently. However, this optimal performance is an extra layer of complexity that needs to be solved manually and AD takes care of that automatically.

The remaining object to be computed from Eq. (\ref{analytic_der}) is $\gamma$. Since, $A$ is a very large sparse matrix the best way to do determine it is by solving the adjoint system equivalent to $\gamma=\tilde{M}A^{-1}$ given as
\begin{align}\label{adjoint_system}
	A^T\gamma = \tilde{M}^T \text{ with } \gamma\in \mathbb{R}^{N+(L-1)\times L}.
\end{align}

Since $A$ depends on the conductivity $\sigma$ this system needs to be solve once at each iteration of the inverse solver.

Finally, a formula for the derivatives in Eq. (\ref{jacobian}) is obtained after solving the adjoint system (\ref{adjoint_system}) and computing the derivative of $B^1$ as in (\ref{b1_der_w}). The derivatives are compactly given through the formula
\begin{align}
	\frac{\partial V_n}{\partial w} = - \gamma^T \begin{bmatrix}
		\frac{\partial B^1}{\partial w} & 0 \\[2pt] 0 & 0 
	\end{bmatrix}
	\theta_n.
\end{align}

Through this demonstration, we have seen that it can be very tedious to deduce and implement the analytical derivatives for complex problems, like ours. For simple functions, an analytical derivative in compact form takes the lead in efficiency, however we want to experiment with the case of more complex functions.

\subsection{Automatic differentiation method}

In this section, we introduce how to apply JAX automatic differentiation toolbox \cite{jax} to obtain the Jacobian. Further details about the inner workings of AD and JAX are explained in Appendix C. 

Since our direct operator $\Sim$ has more output variables than input variables we note that the most-efficient AD mode is the forward-mode.

The implementation of a differentiable simulator $\Sim$ means we can simply use JAX AD to compute the derivatives. Our $\Sim$ operator is differentiable with respect to the parameterization variables $(r,c_x,c_y, \sigma_{in}, \sigma_{out})$ that define the anomaly, as introduced in section 3.2.

This preparation are a requirement for both derivative methods, but now the derivative computation with AD is simply implemented through JAX.

To do so, we implement a routine that defines the direct operator $\Sim$ given in Eq. (\ref{para_Sim}). The implementation is established through the solution of the direct problem through FEM, that we here hide as the \textit{simulator} method. Listing \ref{list:direct_code} provides the routine with all of these in mind.
\begin{listing}[!ht]
	\begin{minted}{Python}
import jax
	
def direct_operator(anomaly_parameters):
 """Simulate measurements for given
 input function with JAX.

 Args:
 anomaly_parameters: Array of shape
  (5,) with parametrization variables
  of circular anomalies.

 Returns:
 measurements: Array of shape
  (nmb_electrodes(nmb_electrodes-1),) that
  contains the voltage measurements for all
  current patterns.
 """

 # Compute measurements
 measurements = simulator(anomaly_parameters)

 return measurements
\end{minted}
\caption{Definition of the direct operator through a general simulator method.}
\label{list:direct_code}
\end{listing}

In order to compute the Jacobian defined in Eq. (\ref{jacobian}) with JAX  one only needs to call \texttt{jax.jacfwd(direct\_operator)} for our direct operator as in Listing \ref{list:jax_jacobian}.
\begin{listing}[!ht]
\begin{minted}{Python}  
def jacobian(anomaly_parameters):
  """Compute Jacobian with JAX AD
  
  Args:
  anomaly_parameters: 1d array of shape
   (5,) with parametrization variables
   of circular anomalies.
  		
  Returns:
  Jacobian matrix of shape
  (nmb_electrodes(nmb_electrodes-1), 5).
  """

  # Define the jacobian through forward-mode
  jacobian = jax.jacfwd(direct_operator)

  return jacobian(anomaly_parameters)
\end{minted}
\caption{Computation of the Jacobian matrix through JAX automatic differentiation toolbox.}
\label{list:jax_jacobian}
\end{listing}

To establish the inverse solver these function definitions are redundant and we can immediately call \texttt{simulator} and \texttt{jax.jacfwd(direct\_operator)} in the inverse solver routine.  This definition is just for visualization purposes in this section.

\section{Experimental setup}

To compare both analytic and automatic differentiation methods, we explore their evaluation at different conductivities, and how they fit in to solve the inverse problem. For the latter, we consider two particular cases for the inverse problem. The first case, that we label as the case of \textbf{fixed conductivities} is simpler. We want to determine only the location parameters $(r, c_x, c_y)$ and we assume the conductivity values inside $\sigma_{\mathrm{in}}$ and outside $\sigma_{\mathrm{out}}$ are fixed. This scenario can represent breast cancer, for example, where we know \textit{a priori} conductivity values of different tissues, and we are only concerned in determining the anomaly location. 

The second case, that we label as the case of \textbf{general conductivities}, we want to determine all parameters $(r, c_x, c_y, \sigma_{in}, \sigma_{out})$. This is a more general scenario where we only know there is a circular anomaly and want to characterize it in terms of location, radius and conductivity.

Recall that we fix a voltage measurement setup to simplify the comparison. Our only interest is to show that AD is as good as analytical methods in terms of solution accuracy. Further, we show that the memory requirements for AD scale reasonably well with the mesh resolution, to show that AD can be effectively implemented in more realistic cases involving more complex scenarios and 3D meshes.

All of the experiments have been run in a machine with the following hardware specifications:
\begin{itemize}
	\item CPU \textit{Intel Core i5-12400F} (released in Q1 2022, 12th gen., 4.4 GHz, 6 cores, 12 threads, 64 GB RAM);
	\item GPU \textit{NVIDIA GeForce RTX 3070} (released in Q4 2020, 6144 CUDA cores, 8 GB memory).
\end{itemize} 

We chose this machine because it has typical med-range specs and can be considered as a good example of an affordable solution for the numerical computation, compatible with the lower cost of EIT.
We remark that besides automatic differentiation, JAX excels in optimizing the performance for a given hardware. Therefore, we have not performed any specific optimization, but appropriate care as been taken throughout the implementation.

\subsection{Establish a ground truth}

In order to have a ``lab" setup, \textit{i.e.}, one we can control the experiment from start to finish, we define a voltage measurements dataset through simulation. For such, we randomly initialize our conductivity parameterization under a certain range of parameters and determine their respective voltage measurements $m$.

To test new inverse solvers we need to generate measurements with the highest resolution possible to avoid the so-called \textit{inverse crimes}. Such crimes occur by using the same resolution to obtain $m$ and $\Sim$ operator computationally. By doing it, we do not account for errors arising from the approximate nature of the direct solver, which occurs when using true measurements obtained by a real-world measuring device, which adequately we can think as having infinite resolution. As such, we need to choose a higher mesh resolution for $m$ than for $\Sim$ operator, since they are obtained both through FEM.

With this in mind we generate our ground truth dataset of voltage measurements with the highest possible resolution for our hardware specifications. In our work, it was established with a FEM mesh of $5815$ elements that is set accordingly to have each element with a edge length of $h=0.035$ relative to the domain size.

Furthermore, we generate the dataset through the following random initialization of the anomaly parameters:
\begin{itemize}
	\item Uniformly generate conductivity centers anywhere inside the disk domain $\Omega=B_1(0)$ with radius 1. Hence, we use polar coordinates to generate the centers. To start we  uniformly generate an angle between $[0, 2\pi]$. Then, we uniformly generate a value in $[0, 1]$ to obtain a radius sample by taking square root of it. Joining both through polar coordinates gives an almost uniformly sampled set of 2D points inside $\Omega$;
	\item Uniformly generate an anomaly radius, taking into consideration the center position generated on the previous point, so that anomalies are strictly in $\Omega$. As such, for each center we select the anomaly radius uniformly from $[0.1, 1-|c|]$, where $|c|$ is the distance from center to origin;
	\item Uniformly generate conductivity values inside $\sigma_{\mathrm{in}}$ from $[1, 1.6]$ S/m and outside $\sigma_{\mathrm{out}}$ from $[0.6, 1.]$ S/m. Such values do not encapsulate any particular medical or industrial scenario.
\end{itemize}

Our model assumes that contact impedances on each electrode are fixed and have value $z=5\times 10^{-6} \Omega\cdot$.

In fact, we generate two separate datasets each with $1000$ cases. One for the case of \textbf{fixed conductivities} where we randomly generate $1000$ anomalies and compute the respective measurements with fixed conductivity value inside of $\sigma_{\mathrm{in}}=1.4$ S/m and outside of $\sigma_{\mathrm{out}}=0.7$ S/m. Another for the case of \textbf{general conductivities} where we randomly generate $1000$ anomalies and compute their measurements as described above.

Furthermore, we provide an initial sanity check for the general dataset. We verified that the Jacobian computed through both methods matches with minimal error margin, which may arise due to round off errors. This analysis is presented in Appendix D.

\section{Results} 

In order to solve the inverse problem for the two cases described above, we use a FEM mesh with $5210$ elements set by $h=0.037$ to define the $\Sim$ operator, in order to avoid inverse crimes. Our chosen inverse solver is the Levenberg-Marquardt method with a line search algorithm on each iteration. Further, we establish two stopping criteria based on a maximum number of iterations equal to $20$ and a relative mean squared loss \begin{align}\frac{1}{2}\frac{\left\|\Sim(\sigma)-m^{\text{true}}\right\|_2^2}{\|m^{\text{true}}\|_2^2}< \xi
\end{align} 
with a feasible threshold of $\xi = 0.001$. This choice was established empirically, since after that it becomes hard to improve the anomaly reconstruction.

Let $\sigma_{\text{AD}}$ and $\sigma_{\text{AN}}$ be the solutions obtained through the inverse solver with the different methods to compute the derivative. In order to verify the effectiveness of AD in solving the EIT inverse problem we evaluate how $\sigma^{\text{AD}}$ and $\sigma^{\text{AN}}$ compare with the true solution $\sigma^{\text{true}}$ and how they compare with each other. This evaluation is based on the mean squared error between the anomalies, \textit{i.e.}, for two different anomaly parameterizations $\sigma_1, \sigma_2$ we evaluate $$\text{MSE}(\sigma_1,\, \sigma_2):=\|\sigma_1-\sigma_2\|_2.$$ In essence, we compute
$\text{MSE}(\sigma^{\text{true}},\,\sigma^{\text{AD}}),$ $ \text{MSE}(\sigma^{\text{true}}, \sigma^{\text{AN}}),$ $  \text{MSE}(\sigma^{\text{AD}}, \,\sigma^{\text{AN}}).
$
Then, we perform an analysis of the mean squared errors by computing simple statistics of the mean, variance, maximum and minimum error, and by plotting the histogram with a logarithmic scale in the x-axis. 

We remark that the following analysis is focused on a general analysis on the reconstructions obtained through the different methods and does not verifies the nature of the errors obtained, \textit{i.e.}, we do not check if the errors are occurring for one specific parameter or for small/large values of those same parameters.  

\subsection{Case 1: Fixed Conductivities}

In this case our goal is to determine the anomaly parameterized by $\sigma^{\text{true}}=(r, c_x, c_y)$, since we know \textit{a priori} that the conductivity inside and outside are $\sigma_{in}=1.4$ S/m and $\sigma_{out}=0.7$ S/m, respectively. Here, we denote $\sigma^{\text{true}}$ as the conductivity we aim to discover and $m^{\text{true}}$ for the respective measurements.

We start from our measurements dataset for the \textbf{fixed conductivities} with the set of $1000$ voltage measurements corresponding to different anomalies. This number of experiments was constrained by time and hardware capabilities.  

\begin{figure}[!ht]
	\includegraphics[scale=0.4]{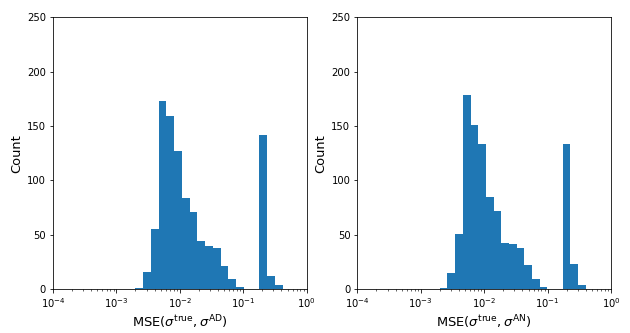}
	\caption{Histogram of the mean squared errors of fixed conductivities, case 1, comparing the reconstructed anomalies obtained through the different derivative methods with the true anomalies.}
	\label{fig:hist_mse_fixed_conductivities}
\end{figure}

\begin{figure}[!ht]
	\centering
	\includegraphics[scale=0.4]{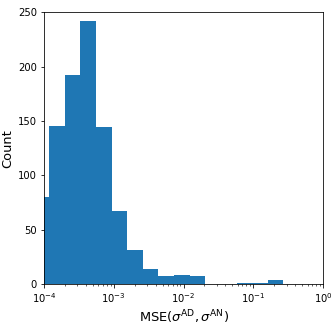}
	\caption{Histogram of the mean squared errors of fixed conductivities, case 1, comparing the reconstructed anomalies.}
	\label{fig:hist_mse_fixed_conductivities_}
\end{figure}

The statistical analysis for this case is given in Table \ref{table:statistic_analysis_fixed_conductivities} and the histogram for the different mean squared errors are in Fig. \ref{fig:hist_mse_fixed_conductivities} and \ref{fig:hist_mse_fixed_conductivities_}.

\begin{table}[!ht]
	\centering
    \vspace{0.25cm}
	\begin{tabular}{c | c  c  c  c}
		& Mean & $S^2$ & Max. & Min. \\
		\hline
		$\text{MSE}(\sigma^{\text{true}}, \sigma^{\text{AD}})$ & 0.0456 & 0.0059 & 0.4177& 0.0020\\
		$\text{MSE}(\sigma^{\text{true}}, \sigma^{\text{AN}})$ & 0.0455 & 0.0057 & 0.4007 & 0.0020 \\
		$\text{MSE}(\sigma^{\text{AD}}, \sigma^{\text{AN}})$ &  0.002 & 2.64e-4& 0.2702 & 1.51e-5 \\
		\hline
	\end{tabular}
	\caption{Statistics of mean squared errors of fixed conductivities, case 1, that compares the reconstructed conductivities obtained through the different derivative methods with the true anomalies.}
	\label{table:statistic_analysis_fixed_conductivities}
\end{table}

The histogram presented in the Fig. \ref{fig:hist_mse_fixed_conductivities} shows that the distribution of the mean squared errors $\text{MSE}(\sigma^{\text{true}}, \sigma^{\text{AD}})$ and $\text{MSE}(\sigma^{\text{true}}, \sigma^{\text{AN}})$ is similar. Notice that the mean squared errors in both cases are concentrated around $10^{-2}$ with a set of outliers with error higher than $0.1$. However, this outliers occur in the same proportion for both methods.  In analysis, this shows that the inverse solver with automatic differentiation matches that with the analytic derivative.

In Fig.~\ref{fig:hist_mse_fixed_conductivities_} the histogram presents the distribution of the mean squared errors between reconstruction $\text{MSE}(\sigma^{\text{AD}},\sigma^{\text{AN}})$ and one can see that it is highly concentrated around $10^{-3}$. There are some different reconstructions between the methods, but their error is in the order of $0.1$. Again, this highlights the effectiveness of AD compared with the analytic method. However, there are some outliers that shows divergence in the reconstructions between both methods. These errors seem to be related with round-off errors when we combine this analysis with the sanity check for the Jacobian.

To complete the discussion of this case, we allude to the statistics in Table \ref{table:statistic_analysis_fixed_conductivities}. We point to the mean and variance of the different mean squared errors. This shows that on average the reconstruction obtained with AD is much closer with the analytic one than with the true anomalies. Furthermore, the variance between these reconstructions is very small. Once again it shows the effectiveness of AD to match the analytic derivative method and that other inverse solver methods need to be improved in order to obtain better reconstruction results.

\subsection{Case 2: General Conductivities}

For this case the objective is to determine the general anomaly parameterization given by $\sigma^{\text{true}}=(r, c_x, c_y, \sigma_{\mathrm{in}}, \sigma_{\mathrm{out}})$. Again, we denote $\sigma^{\text{true}}$ as the conductivity we aim to discover and $m^{\text{true}}$ for the respective measurements.

We start from the measurements dataset for the \textbf{general conductivities} with the set of $1000$ voltage measurements corresponding to the different anomalies. Recall, that in this generation we have assumed that $\sigma_{\mathrm{in}}$ is always greater than $\sigma_{\mathrm{out}}$.

\begin{figure}[!ht]
	\centering
	\includegraphics[scale=0.4]{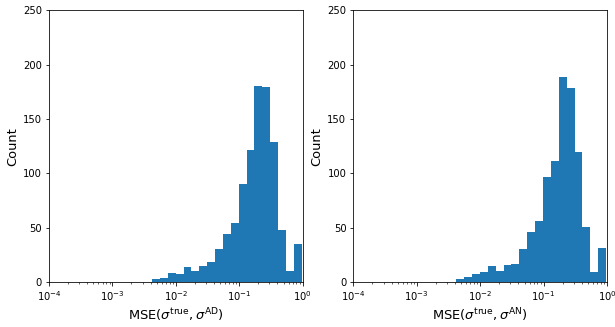}
	\caption{Histogram of the mean squared errors of general conductivities, case 2, that compares the reconstructed anomalies obtained through the different derivative methods with the true anomalies.}
	\label{fig:hist_mse_general_conductivities}
\end{figure}

\begin{figure}[!ht]
	\centering
	\includegraphics[scale=0.4]{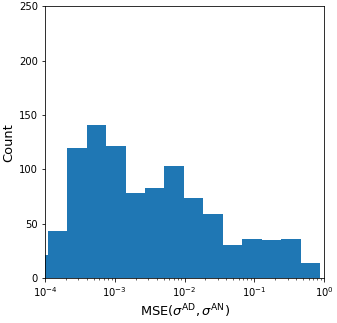}
	\caption{Histogram of the mean squared errors of general conductivities, case 2, comparing the reconstructed anomalies.}
	\label{fig:hist_mse_general_conductivities_}
\end{figure}

The statistical analysis for this case is given in Table \ref{table:statistic_analysis_general_conductivities} and the histogram for the different mean squared errors are in Figs. \ref{fig:hist_mse_general_conductivities} and \ref{fig:hist_mse_general_conductivities_}.

\begin{table}[!ht]
	\centering
	\vspace{0.25cm}
	\begin{tabular}{c | c c c c}
		& Mean & $S^2$ & Max. & Min. \\
		\hline
		$\text{MSE}(\sigma^{\text{true}}, \sigma^{\text{AD}})$ & 0.2264 & 0.0292 & 0.9698& 0.0042 \\
		$\text{MSE}(\sigma^{\text{true}}, \sigma^{\text{AN}})$ & 0.2215 & 0.0273& 0.9706 & 0.0042 \\
		$\text{MSE}(\sigma^{\text{AD}}, \sigma^{\text{AN}})$ &  0.039 & 0.0134& 0.8838 & 4.4e-6 \\
		\hline
	\end{tabular}
	\caption{Statistics of mean squared errors of general conductivities, case 2, that compares the reconstructed conductivities obtained through the different derivative methods with the true anomalies.}
	\label{table:statistic_analysis_general_conductivities}
\end{table}

The histogram presented in the Fig. \ref{fig:hist_mse_general_conductivities} shows that the distribution of the mean squared errors $\text{MSE}(\sigma^{\text{true}}, \sigma^{\text{AD}})$ and $\text{MSE}(\sigma^{\text{true}}, \sigma^{\text{AN}})$ is similar. In analysis, this shows that the inverse solver with automatic differentiation matches that with the analytic derivative. Further, notice that the mean squared errors in both cases are concentrated around $10^{-1}$. In fact by setting a threshold, we verified that there are at most $50$ reconstructions for both methods where the mean squared error with the true anomaly is higher than $0.5$, which together with the histograms shows that the vast majority of reconstructions is successful.

Furthermore, the histogram in Fig. \ref{fig:hist_mse_general_conductivities_} that presents the histogram of  $\text{MSE}(\sigma^{\text{AD}},\sigma^{\text{AN}})$ shows that the errors between reconstructions are more concentrated around the interval $[10^{-4}, 10^{-2}]$. Again, this highlights the equivalence of AD compared with the analytic method. However, there are some outliers that shows divergence in the reconstructions between both methods. Combining this analysis with the sanity check for the Jacobian it reveals that this might occur due to round-off errors.

To complete the discussion of this case, we allude to the statistics Table \ref{table:statistic_analysis_general_conductivities}. The only aspect we would like to point out here is the mean of the different mean squared errors. This shows that on average the reconstruction obtained with AD is much closer with the analytic one than with the true anomalies. Once again it shows the effectiveness of AD to match the analytic derivative method and that other inverse solver methods need to be improved in order to obtain better reconstruction results. 

\subsection{Computational performance of AD}

The viability of AD also depends of its scaling capabilities. Namely, we want to understand if increasing the number of mesh elements, and therefore the resolution and accuracy of the FEM turns AD unfeasible. This is relevant because AD requires the construction of a computational graph for the direct problem and then applies the chain-rule throughout the nodes of the graph to compute the derivatives. As the number of mesh elements increases the computational graph becomes larger and can be unfeasible to use for it to compute the derivatives.

\begin{figure}[!ht]
	\includegraphics[scale=0.35]{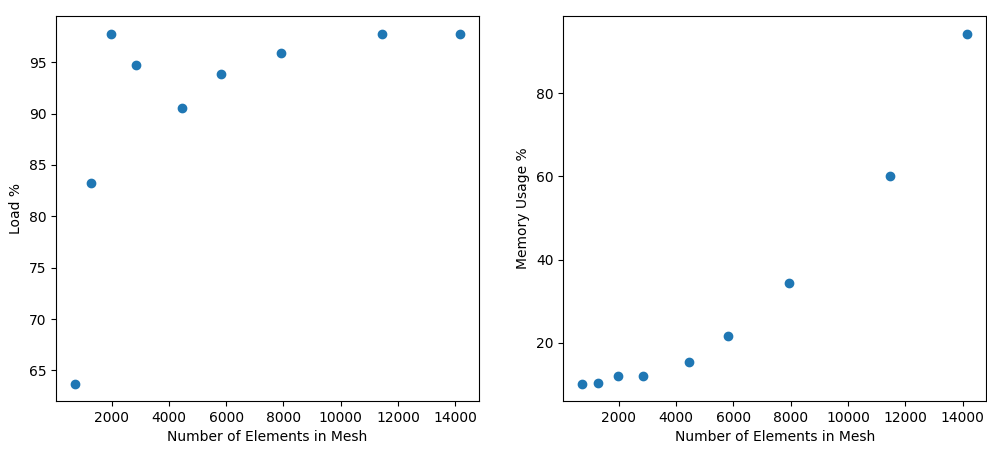}
	\caption{Percentage of GPU load and memory usage with respect to the number of mesh elements.}
	\label{fig:gpu_memory_load}
\end{figure}

\begin{figure}[!ht]
	\centering
	\includegraphics[scale=0.35]{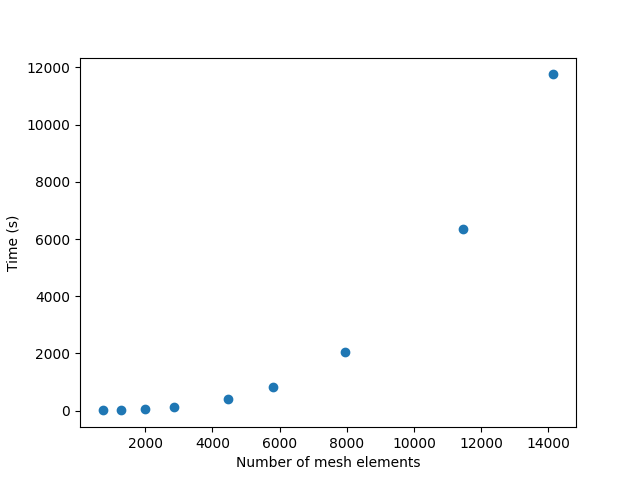}
	\caption{Time (s) elapsed to compute Jacobian matrices for 100 random anomalies with respect to the number of mesh elements.}
	\label{fig:gpu_time_}
\end{figure}

In order to understand this behavior, we compute for ten different mesh sizes the Jacobian for $100$ distinct general anomalies, randomly generated as described before. For each mesh size we measure the average GPU memory and load usage through the Python package \texttt{GPUtil}. In Fig. \ref{fig:gpu_memory_load} we plot the average of GPU load and memory usage percent for each of the different mesh resolutions and in Fig. \ref{fig:gpu_time_} we plot the time that took to compute the Jacobian matrices with respect to each mesh resolution.

It is clear from both figures the growth in GPU memory usage and time to execute this experiment. Moreover, for meshes with more than 15000 elements we require more than 8Gb of GPU memory. As of now, we cannot understand the order of growth and further experiments with finer resolution are needed. 

\section{Conclusion}

In this paper we have compared the effectiveness of AD to solve inverse problems against classical methods with analytical formulations of the derivative. We have shown how to adequately construct a FEM differentiable simulator in the context of inverse problems. We successfully introduced automatic differentiation for solving inverse problems in an optimization framework, in particular, Electrical Impedance Tomography. We have shown that AD provides a simple way of computing derivatives of complex operators, for example, arising from solutions of PDEs, with respect to a set of parameters.

We have shown that AD is indeed effective to compute the derivatives, since it matches the analytical computation up to minimal error. Further, it was used to solve the Electrical Impedance Tomography inverse problem and we shown that it is even superior to analytical methods, in terms of time and resources. 

The analytical formulation is nothing more than an application of differentiation rules to the FEM formulation of the direct operator. By construction AD essentially executes the same process, but automatically. As such, AD and the analytical formulation can be even performing the same operations, but the fact that AD is a plug-and-play tool makes it advantageous to use for complex operators. 

Moreover, it has proven more efficient since it takes less time on average to solve any particular EIT problem, when compared with the analytical formulation in our case study and scales well with the mesh resolution. This indicates that with the right hardware AD can be efficiently executed for large-scale problems.

With this tool, we can cast our focus into an efficient implementation of the direct problem solvers, which is way more understood in literature, and on the methods to solve the inverse problem. It allows freedom to experiment and deal with difficult equations, without much thought, bringing focus to the practical application at hand.

Further, we expect that AD extends nicely to higher dimensions, while the analytic formulation will require some re-implementation to accommodate the three dimensional shapes of anomalies.

Future studies are interested in testing how AD easily handles different shapes of anomalies, as well as 3D models.

\begin{ack}
The work of I. Pombo was supported by FCT through CIDMA and projects UIDB/04106/2020, UIDP/04106/2020 and the PhD Scholarship SFRH/BD/143523/2019. 

This work was developed during a research internship at Inductiva Research Labs from March 2022 to Jan 2023. The first author would like to thank the entire Inductiva team for the continuous support and encouragement  provided during the entire period of the internship and in particular thank Hugo Penedones, Fábio Cruz, David Lima and David Carvalho for their comments and constructive feedback given over the several versions of this manuscript.
\end{ack}

\appendix

\section{FEM formulation of the Direct Problem}

In this appendix, we take a deeper dive into the Finite Element Method applied to Electrical Impedance Tomography. In this paper, we have used the Complete Electrode Model (CEM) \cite{cheng} to understand how current propagates inside the body. Starting from its formulation, that we have introduced in section 2.2. and recall here, we derive its weak formulation and apply FEM to obtain a system of linear equations.

The CEM takes into account the finite nature of electrodes, current injection through them and electro-chemical effects happening between skin and electrode surface. 

Let $\Omega$ describe the subject region we are evaluating. To establish a measurement setup, we attach $L$ electrodes at the subject boundary $\partial\Omega$. Through them we apply an electrical \textit{current pattern} $I=\left(I_1,...,I_L\right)$ into $\Omega$. The objective is to find the electrical potential $u$ inside and the voltages at electrodes $V=\left(V_1,...,V_L\right)$
that fulfill the system of equations describing the Complete Electrode Model:
\begin{align}\label{CEM_}
\begin{cases}
\nabla\cdot(\sigma\nabla u)=0, &\quad\text{ in } \Omega,\\
\int_{E_l} \sigma\frac{\partial u}{\partial \nu}\, dS = I_l, &\quad l=1,2,...,L \\
\sigma\frac{\partial u}{\partial\nu} = 0,&\quad \text{ in } \partial\Omega\setminus\cup_{l=1}^L E_l \\
\left. u + z_l\sigma\frac{\partial u}{\partial\nu}\right|_{E_l} = V_l, &\quad l=1,2,...,L
\end{cases}
\end{align}
where $\sigma$ is the conductivity distribution. 

The first equation represents electrical \textit{current diffusion}. The second and third define the insertion of current through electrodes, meaning that current spreads through the whole electrode before being inserted into the domain and in regions without electrode there isn't current flowing. Finally, the last equations models the electrochemical effects at interface of skin-electrode, with $z_l$ designated as \textit{contact impedance} represent the resistance at that interface.

To ensure the existence and uniqueness of a solution, the current pattern must satisfy Kirchoff's law and we fix a \textit{reference voltage condition}:
\begin{align}\label{uniq_cond_}
\sum_{l=1}^L I_l=0, \quad \text{ and }\quad \sum_{l=1}^L V_l=0.
\end{align}

In order to apply  Finite element method, we introduce the variational equation that describes fully (\ref{CEM_}). In \cite{somersalo} it has been derived and shown that $(u, V)$ is a weak-solution of (\ref{CEM_}) if for all $(w, W)\in H^1(\Omega)\times \mathbb{R}^L$ we have:
\begin{align}\label{variational_eq}
\int_{\Omega} &\sigma\nabla u\cdot\nabla vdx + \sum_{l=1}^L \frac{1}{z_l}\int_{E_l} \left(u-V_l\right)\left(w-W_l\right) dS\nonumber \\ 
&\hspace{5cm} = \sum_{l=1}^L I_lW_l
\end{align}

This formulation joins every condition of (\ref{CEM_}) together into one equation, which allows the simplification into a linear system of equations. The first integral describes the propagation of current throughout the domain, while the second represents skin-electrode interface condition and the right-hand side explains the insertion of current.

\subsection{FEM for Complete Electrode Model}

FEM allows transforming the continuous problem, described by the variational equation (\ref{variational_eq}), into a discrete system of equations that can be handled by linear algebra methods. 

A detailed explanation is provided in any FEM book, and specifically for EIT \cite{siltanen}. Further, in Appendix we provide an explanation of each step for complete understanding of those interested. Here, we only briefly describe some of the parts required for our exposition.

In this section, we briefly describe how to apply FEM in EIT. First, we  remark that many variants can arise due to possible different assumptions made.

First, we discretized the subject domain $\Omega$ into smaller elements. Next, we approximate our solutions $u, U$ by
\begin{align}\label{approx_sol_}
u^h(x, y) &= \sum_{i=1}^N \alpha_i\phi_i(x, y) \\
 V^h &= \sum_{k=1}^{L-1} \beta_k\eta_k,
\end{align} 
where $\phi_i, \eta_k$ are basis functions. In particular, the $\eta_k$ are defined through $\eta_1=(1,-1,0,...,0)^T,\, \eta_2=(1,0,-1,0,...,0)^T, \,..., \eta_{L-1}=(1,0,...,0,-1)^T\in \mathbb{R}^{L}.$ This choice ensures that reference voltage condition (\ref{uniq_cond_}) is fulfilled. Further, $N$ in (\ref{approx_sol_}) corresponds to the number nodes forming the finite element mesh.

The approximate solutions $u^h$ and $V^h$ for the direct problem are fully determined by the coefficients
\begin{align}
&\alpha=[\alpha_1,..., \alpha_N] \in \mathbb{R}^N,\\ &\beta=[\beta_1,...,\beta_{L-1}]\in\mathbb{R}^{L-1}.
\end{align}

FEM allows us to obtain a system of linear equations characterizing them. This is achieved by inserting $(u^h, V^h)$ into the variational equation (\ref{variational_eq}), together with different choices of $(v, V)=(\phi_i, \eta_j)$. Gathering all possibilities leads to a linear system of equations:
\begin{align}\label{system_eq_}
A\theta = \tilde{I},
\end{align} 
where $\theta = [\alpha,\beta] \in \mathbb{R}^{N+L-1}$ and $\tilde{I}$ is described through the current pattern $I$ applied at the electrodes as follows:
\begin{align}
\tilde{I}=\left[\overrightarrow{0}, I_1-I_2, I_1-I_3,...,I_1-I_L\right]\in\mathbb{R}^{N+(L-1)}.
\end{align}

The stiffness matrix $A$ can be computed in terms of four blocks:
\begin{align}\label{stiffness_matrix_}
A=\begin{pmatrix}
B^1+B^2 & C \\
C^T & D
\end{pmatrix}.
\end{align}

Each term is defined through integration over the domain and over the electrodes like:
\begin{align}\label{matrices_def}
B^1_{ij}&=\int_{\Omega} \sigma\nabla\phi_i\cdot \nabla\phi_j\, dx,\quad i, j = 1,2,..., N \\[2pt]
B^2_{ij} &= \sum_{l=1}^{L} \frac{1}{z_l}\int_{E_l} \phi_i\phi_j\, dS, \quad i, j = 1,2,..., N \\[2pt]
C_{ij} &= -\left[\frac{1}{z_1} \int_{E_1} \phi_i\, dS - \frac{1}{z_{j+1}}\int_{E_{j+1}} \phi_i\, dS\right], \nonumber\\
&\hspace{2cm}i=1,2,...,N,\; j=1,2,...,L-1\\[2pt]
D_{ij} &= \begin{cases}
\frac{E_1}{z_1},& \quad i\neq j\\
\frac{E_1}{z_1}+\frac{|E_{j+1}}{z_{j+1}},& \quad i=j
\end{cases}, i,j=1,...,L-1,
\end{align}
with $|E_j|$ being the electrode area. 

The derivation of each block arises from application of two different basis functions on the weak formulation. A full description was done in \cite{siltanen}.

After solving the system for $\theta$, the voltages $V^h$ are obtained by multiplication with the basis functions matrix $M$ defined as:
\begin{align}\label{basis_functions_matrix}
M = \begin{bmatrix}
1 & 1 & 1&\hdots & 1 \\
-1 & 0 & 0&\hdots& 0 \\
0 & -1& 0 &\hdots & 0 \\
\vdots & \vdots &\vdots & \ddots & \vdots \\
0 & 0& 0 & \hdots & -1
\end{bmatrix}
\end{align}
through $$V^h=M\beta.$$

One detail we want to point out regarding FEM implementation concerns conductivity parameterization. For computational purposes, we assume that $\sigma$ is piece-wise constant, meaning that at each mesh element is constant, and thus mathematically defined as:
\begin{align}\label{conductivity_discrete_}
\sigma(x, y) = \sum_{k=1}^{K} \sigma_k\chi_k(x, y),
\end{align}
where $K$ is the total number of elements and $\chi_k$ is the indicator function of the $k$-th element.

In this sense, matrix $B^1$ simplifies to
\begin{align}\label{discrete_B1_}
B^1_{ij} = \sum_{\{k:\, i,j\in T_k\}} \sigma_k \int_{T_k} \nabla\phi_i\cdot\nabla\phi_j\, dx
\end{align}

The parameterization of $\sigma$ is essential to compute the voltages variation $V^h$ with respect to a conductivity variation, i.e., the derivative. If a parameterization was not applied to $\sigma$, then it would be described as a function from $\Omega$ to $\mathbb{R}$. For the latter case a derivative still exists, but it is more theoretically described, see \cite{harrach}. 
\\

\subsection{Implementation Details}

For implementation purposes we restrict ourselves to two-dimensions even though the above formulation also holds for further dimensions.

The first implementation decision is about space discretization. For simplicity sake, our choice of mesh generator is \textit{DistMesh algorithm}, developed by Per-Olof Persson and Gilbert Strang \cite{persson}. The elements are triangles and the algorithm has been adapted to consider $L$ equidistant electrodes, with a pre-defined size, at the surface $\partial\Omega$.

Secondly, we need to define our basis functions. We choose piece-wise linear functions, and therefore, for each triangle element any basis function is linearly defined as:
$$\phi_i(x,y)=a_i+b_ix+c_iy.$$

Moreover, the basis function are obtained in correspondence to a mesh node $(x_j, y_j)$ through the condition:
\begin{align}
\phi_i(x_j,y_j) =  \begin{cases}
1,& \quad i=j \\
0,& \quad i\neq j
\end{cases}.
\end{align}

Since, for every other node the function will be zero, it holds that for every triangle that does not have $i$ as a node, $\phi_i\equiv0$ there. This simplifies the computations of all the matrices, since most entries will be $0$, due to non-intersection of most basis functions supports. As such, $A$ is sparse.

Moreover, due to the equation nature being elliptic, the stiffness matrix $A$ is positive-definite. As such, the most appropriate system of equations solver is the \textit{Conjugate Gradient method (CG)}. 

\section{Derivation of Levenberg-Marquardt method}

A simple method for inverse problems under such an optimization framework is \textit{Levenberg-Marquardt method.}

It is a general method since it is independent of the simulator and the method used for differentiating it. As such, it allows us to demonstrate the effectiveness of various methods to compute the derivatives, in particular, of Automatic Differentiation. 

We hereby assume that $\sigma$ is discretely given by a parameterization, i.e., $\sigma\in\mathbb{R}^p$. This simplifies simulation and, more importantly, the derivatives computation process which is now done with respect to each variable $\sigma_i,\, i=1,...,p$. An example is seen in Figure \ref{fig:cancer_id} where $\sigma=(\sigma_1,\sigma_2)$.

The minimization problem is given as
\begin{align}\label{min_2}
\min_{\sigma} \frac{1}{2}\left\|\Sim(\sigma) - m^{\text{true}}\right\|_2^2,
\end{align}
where $m^{\text{true}}$ is a set of true measured voltages with respect to $N$ currents applied, as already introduced. 

We iteratively improve an approximate solution of the minimization problem (\ref{min}) through 
\begin{align}
\sigma^{k+1}=\sigma^k + \delta\sigma^{k}
\end{align}
where $\delta\sigma^k$ is an \textit{update step} and $\sigma^k$ is the current approximate solution. This process is done until a satisfactory solution is found.

Each method to solve the minimization problem is defined by the update step $\delta\sigma^k$ computation.

Levenberg-Marquardt is a second order quasi-newton method, that approximates the Hessian through an identity regularization. In this sense, the update rule is given as follows:
\begin{align}\label{LM_update_rule}
\resizebox{0.9\hsize}{!}{$
	\delta\sigma^{LM}=-\left[J(\sigma)^TJ(\sigma)+\lambda_{LM}I\right]^{-1} J(\sigma)^T\left(\Sim(\sigma) - m^{\text{true}}\right)$}.
\end{align}

Here, $J(\sigma)$ denotes the Jacobian of $\Sim$, i.e., a matrix of voltage derivatives with respect to each parameter $\sigma_i$. Further, $\lambda_{LM}$ is a parameter used to approximate the Hessian and that allows for improving the condition number of $J(\sigma)^TJ(\sigma)$. We determine it empirically.

The update rule derivation is given in appendix.

The Levenberg-Marquardt method is a particular type of quasi-Newton methods. We start by deducing the general form of quasi-Newton methods and there after funnel on our chosen method.

We hereby assume that $\sigma$ is discretely given by a parameterization, i.e., $\sigma\in\mathbb{R}^p$. This simplifies simulation and, more importantly, the derivatives computation process which is now done with respect to each variable $\sigma_i,\, i=1,...,p$. An example is seen in Figure \ref{fig:cancer_id} where $\sigma=(\sigma_1,\sigma_2)$.

The minimization problem is given as
\begin{align}\label{min_2_}
\min_{\sigma} \frac{1}{2}\left\|\Sim(\sigma) - m^{\text{true}}\right\|_2^2,
\end{align}
where $m^{\text{true}}$ is a set of true measured voltages with respect to $N$ currents applied, as already introduced. 

Denote by $\mathcal{L}(\sigma)$ the loss function in (\ref{min_2_}). Then, assuming that we have an initial guess $\sigma_0$, we can re-write (\ref{min_2}) as
\begin{align}\label{loss_}
\mathcal{L}(\sigma+\delta\sigma) = \frac{1}{2}\left\|\Sim(\sigma_0+\delta\sigma)-m^{\text{true}}\right\|_2^2
\end{align}
with an intent to minimize with respect to the parameter variation $\delta\sigma$, which we designate by \textit{update step}. Thereafter, applying this iteratively we approximate our solution through
\begin{align}
\sigma_{k+1}=\sigma_k + \delta\sigma_{k}
\end{align}
until a satisfactory solution is found.

The Levenberg-Marquardt Algorithm is essential for to compute the update step $\delta\sigma$. Taylor expansion of (\ref{loss_}) up to quadratic term is given by
\begin{align}\label{L_sigma}
\mathcal{L}(\sigma+\delta\sigma)= \mathcal{L}(\sigma) + \mathcal{L}'(\sigma)\delta\sigma + \frac{1}{2}\mathcal{L}''(\sigma)(\delta\sigma)^2,
\end{align}
where $\mathcal{L}'(\sigma)$ and $\mathcal{L}''(\sigma)$ denotes the gradient and Hessian of the objective function $\mathcal{L}$, with respect to parameters defining $\sigma$.

A minimum with respect to $\delta\sigma$ has gradient zero. Thus, we apply gradient to (\ref{L_sigma})
\begin{align*}
	\frac{\partial \mathcal{L}}{\partial\delta\sigma}(\sigma+\delta\sigma) = \mathcal{L}'(\sigma) + \mathcal{L}''(\sigma)\delta\sigma,
\end{align*}
and setting the gradient equal to zero yields
\begin{align*}
0 &= \mathcal{L}'(\sigma) + \mathcal{L}''(\sigma)\delta\sigma \\
\Leftrightarrow
\delta\sigma &= -\left[\mathcal{L}''(\sigma)\right]^{-1}\mathcal{L}'(\sigma).
\end{align*}

Since only $\Sim$ depends on conductivity parameterization we can compute the gradient and Hessian through:
\begin{align}\label{gradient}
\mathcal{L}'(\sigma)&=J(\sigma)^T \left(\Sim(\sigma)-m^{\text{true}}\right) \\
\label{hessian}
\mathcal{L}''(\sigma) &= J(\sigma)^T J(\sigma) + \nonumber \\ &\hspace{1cm} \sum_i \left[\Sim_i(\sigma)\right]''\left(\Sim_i(\sigma)-m^{\text{true}}_i\right),
\end{align}
where $J$ is the Jacobian of simulated voltages $\Sim(\sigma)$ with respect to the parameterization of $\sigma$.

Up until here the derivation is general for quasi-Newton methods. 

The Levenberg-Marquardt method distinguishes itself from other quasi-Newton methods by avoiding computation second order derivative, substituting it by a scaled identity matrix $\lambda_{LM}I,\,\lambda\in\mathbb{R}^+$, which acts as a regularizer by improving the condition number of the Hessian matrix to be inverted. Now, the update can be computed through:
\begin{align}\label{update_rule}
\resizebox{0.9\hsize}{!}{$
	\delta\sigma_{LM}=-\left[J(\sigma)^TJ(\sigma)+\lambda_{LM}I\right]^{-1} J(\sigma)^T\left(\Sim(\sigma) - m^{\text{true}}\right)$}.
\end{align}  

\section{Automatic Differentiation}

AD is a set of techniques to evaluate the derivative of a function specified by a computer program. No matter how complicated they are, any computer program is based on a simple set of arithmetic operations and functions, like addition, multiplication, trigonometric functions, exponentials, etc. We can encode the derivative rule for all of these simple operations and build up the full derivative of our complex program through the chain-rule. AD evaluates derivatives with exact precision.  

There are two modes for AD implementation: forward-mode and reverse-mode. In any case, they are not hard to implement through operator overloading techniques. The difficult part is to provide an efficient and optimal computation of these modes. However, at the present moment there are great libraries that provide efficient implementations of AD for both modes, like JAX for Python. 

The first step in AD is the creation of a computational graph of our program, that explains the decomposition into simpler operations for which we know the derivative. Let's exemplify for the following function $f(x_1,x_2)=\sin(x_1\cdot x_2) + e^{x_1}$. The first step is to break things apart into the simpler operations:
\begin{align*}
&w_1=x_1, w_2=x_2 \\
&w_3 = w_1\cdot w_2 \\
&w_4 = \sin(w_3) \\
&w_5 = e^{w_1} \\
&w_6 = w_4+w_5 =: f(w_1,w_2)	
\end{align*}

\begin{figure}[!ht]
	\centering
	\includegraphics[scale=0.3]{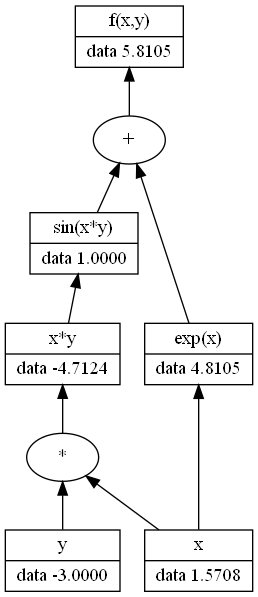}
	\caption{Computational Graph of $f(x_1,x_2)=\sin(x_1 \cdot x_2)+e^{x_1}$ evaluated at $(\pi/2, -3)$.}
	\label{fig:Cgraph}
\end{figure}

This decomposition is more easily visualized through the computational graph in Fig. \ref{fig:Cgraph}.

With the computational graph in mind, forward-mode computes derivatives from bottom-to-top, that is from the variables to output. As such, it allows the derivative computation of all outputs with respect to a single variable. It can evaluate the derivative simultaneously with the function, and thus it is proportional to the original code complexity. In this terms, it is more efficient for functions $f:\mathbb{R}^n\rightarrow \mathbb{R}^m$ with $m>>n$.

Reverse-mode of AD works the other way around, that it is, top-to-bottom. First, it requires a forward evaluation of all the variables, and thereafter it starts computing the derivatives from output values for the variables involved immediately, doing that successively until the input variables. Therefore, it allows evaluation of the gradient of an single output function. As such, it is way more efficient for functions $f:\mathbb{R}^n\rightarrow\mathbb{R}^m$ with $m<<n$. 

A familiar example in these days is neural networks that are described by way more weights that output variables, in this particular example the reverse-mode is known as backpropagation.

One possible limitation to take into account in AD arises from the computational graph we described. Due to the computer program complexity this computational graph can be very expensive to establish and keep in memory. In such scenarios, where the Jacobian is obtained from a very complex graph, instead of a compact formula like analytic formulation, it can take a long time to be evaluated. As such, AD is not a tool to be inserted into play whenever needed and considerations must be made when implementing the $\Sim$ operator, to avoid some of these flaws.

To bypass this problem, JAX can encode loops and conditionals in primitive operations that are inherent from the domain-specific compilers for linear algebra (XLA). Otherwise, the loops are unrolled into a set of operations (may be smaller than the general loop, but) that increases the computational graph size. With the primitives in mind, this will be encoded on the graph with a single operation, for which we already know the derivative.

With AD the focus is completely in an optimal implementation of the $\Sim$ operator, which is essential to obtain a very efficient inverse problem solver (even with analytical computation of derivatives). Thereafter, thinking about both modes, we can apply forward-mode to compute efficiently the derivatives of $\Sim$ with respect to the parameterization $(r,c_x,c_y, \sigma_{in}, \sigma_{out})$. 

Being aware of the inherent problems with both methods is essential for a proper implementation of the inverse solver.

\section{Extended Results}

In this section we present some extra analysis about the Jacobian computation with both methods in order to make a sanity check. 

The sanity check we want to verify is to check if the Jacobian computed through automatic differentiation and the analytic formulation match. This is what we already expect since AD applies the chain-rule of differentiation to FEM, which is exactly what we have done by hand to determine the analytic formulation. The Frobenius norm of the Jacobian difference is given as:
\begin{align*}
	&\left\|J^{AD}-J^{analytic}\right\|_{Fro}, \\
	\text{ where } &\|A\|_{Fro} = \left[\sum_{i, j}^{n, m} |a_{ij}|^2\right]^{1/2}.
\end{align*}

Further, we computed the Jacobian with both methods for $100$ randomly generated general conductivities described in section 3. Thereafter, we compute their difference and applied the Frobenius norm in order to obtain an array with dimension $100$.

\begin{figure*}[!ht]
	\centering
	\includegraphics[scale=0.5]{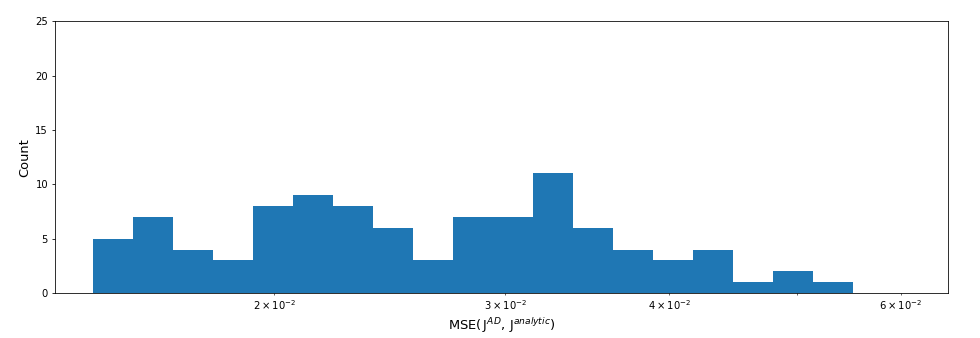}
	\caption{Histogram of Jacobian error with both derivative methods evaluated with Frobenius norm.}
	\label{fig:histogram_jac_loss}
\end{figure*}

To verify the assumption that both should evaluate to almost the same values we make an histogram of the losses and provide some statistics, namely, mean, variance, maximum and minimum. This results are provided in Fig. \ref{fig:histogram_jac_loss} and Table \ref{table:6}.

Statistically we can infer that the Jacobian match closely together with maximum error of $0.0552$ and an average of $0.0271$. Indeed, the histogram confirms that most evaluations are really close together, with only some outliers compared with the overall picture. Further, these outliers might just be rounding off errors and are not worrisome since the error is still considerably small.
\begin{table*}[!ht]
	\centering
	\begin{tabular}{c c c c c}
		& Mean & $S^2$ & Max. Error & Min. Error  \\[0.5ex]
		\hline
		$\left\|J^{AD}-J^{analytic}\right\|_{Fro}$& 0.0271 & 7.94e-05 & 0.0552 & 0.0146
	\end{tabular}
\caption{Statistic analysis of the error between Jacobian matrices obtained through the Frobenius norm.} \vspace{0.25cm}
	\label{table:6}
\end{table*}

\begin{thebibliography}{102}
	
	\bibitem{adler}
	Adler, A., \& Holder, D. (Eds.). (2021). Electrical impedance tomography: methods, history and applications. CRC Press.
	
	\bibitem{jax}
	Bradbury, J., Frostig, R., Hawkins, P., Johnson, M. J., Leary, C., Maclaurin, D., Necula, G., Paszke, A., Vander{P}las, J., Wanderman-{M}ilne, S., \& Zhang, Q. (2018). JAX: composable transformations of Python+ NumPy programs. Version 0.3.13, 5, 14-24.
	
	\bibitem{cheng}
	Cheng, K. S., Isaacson, D., Newell, J. C., \& Gisser, D. G. (1989). Electrode models for electric current computed tomography. IEEE Transactions on Biomedical Engineering, 36(9), 918-924.
	
	\bibitem{tanzer}
	Gençer, N. G., \& Tanzer, I. O. (1999). Forward problem solution of electromagnetic source imaging using a new BEM formulation with high-order elements. Physics in Medicine \& Biology, 44(9), 2275.
	
	\bibitem{harrach}
	Harrach, B. (2021). An introduction to finite element methods for inverse coefficient problems in elliptic PDEs. Jahresbericht der Deutschen Mathematiker-Vereinigung, 123(3), 183-210.
	
	\bibitem{levenberg}
	Levenberg, K. (1944). A method for the solution of certain non-linear problems in least squares. Quarterly of applied mathematics, 2(2), 164-168.
	
	\bibitem{marquardt}
	Marquardt, D. W. (1963). An algorithm for least-squares estimation of nonlinear parameters. Journal of the society for Industrial and Applied Mathematics, 11(2), 431-441.
	
	\bibitem{siltanen}
	Mueller, J. L., \& Siltanen, S. (Eds.). (2012). Linear and nonlinear inverse problems with practical applications. Society for Industrial and Applied Mathematics.
	
    \bibitem{persson}
    Persson, P. O., \& Strang, G. (2004). A simple mesh generator in MATLAB. SIAM review, 46(2), 329-345. 
    
    \bibitem{raissi}
    Raissi, M., Perdikaris, P., \& Karniadakis, G. E. (2019). Physics-informed neural networks: A deep learning framework for solving forward and inverse problems involving nonlinear partial differential equations. Journal of Computational physics, 378, 686-707.   
   
    \bibitem{somersalo}
    Somersalo, E., Cheney, M., \& Isaacson, D. (1992). Existence and uniqueness for electrode models for electric current computed tomography. SIAM Journal on Applied Mathematics, 52(4), 1023-1040.
    
    \bibitem{vauhkonen}
    Vauhkonen, P. J., Vauhkonen, M., Savolainen, T., \& Kaipio, J. P. (1999). Three-dimensional electrical impedance tomography based on the complete electrode model. IEEE Transactions on Biomedical Engineering, 46(9), 1150-1160.
    
    \bibitem{other_vauhkonen}
    Vauhkonen, P. (2004). Image Reconstruction in Three-Dimensional Electrical Impedance Tomography. Kuopio University Publications C. Natural and Environmental Sciences 166.
    
    \bibitem{webster}
    Webster, J. G. (Ed.). (1990). Electrical impedance tomography. CRC Press.
    
    \bibitem{wengert}
    Wengert, R. E. (1964). A simple automatic derivative evaluation program. Communications of the ACM, 7(8), 463-464.
    
\end{thebibliography}
\end{document}